\documentclass[a4paper,10pt]{article}
\usepackage[round]{natbib}
\usepackage[utf8x]{inputenc}
\usepackage{threeparttable}
\pdfoutput=1
\usepackage{verbatim}
\usepackage{epsfig}
\usepackage{epsf}
\usepackage{epstopdf}
\usepackage{graphicx}
\usepackage{textcomp}
 \usepackage{amsmath}
\usepackage{color}
\usepackage{multirow}
\definecolor{Blue}{rgb}{0.3,0.3,0.9}
\definecolor{red}{rgb}{1,0,0}
\usepackage{pgfplots}
\usepackage{tikz}
\usepackage{caption}
\usepackage{subcaption}

\title{\textbf{\fontsize{15}{15}\selectfont Kernel Atomicity}}\author{\smallskip \textbf{\fontsize{11}{11}\selectfont Stephen G. Odaibo$^{1,2}$}\\ \textbf{\fontsize{10}{10}\selectfont M.S.(Math), M.S.(Comp. Sci.), M.D.}\\ \newline\newline \small{\textsl{\fontsize{9}{9}\selectfont $^1$Quantum Lucid Research Laboratories,}}\\  \small{\color{blue}{\fontsize{7}{7}\selectfont stephen.odaibo@qlucid.com}}\\\newline\newline \\\small{\textsl{\fontsize{9}{9}\selectfont $^2$University of Michigan --Ann Arbor}}\\  \small{\textsl{\fontsize{7}{7}\selectfont Ophthalmology and Visual Sciences}} }
\date{}

\date{}

\begin{document}

\maketitle

\begin{abstract}

 Here, we show that the first isomorphism theorem, the orbit-stabilizer theorem, and the non-uniqueness of solutions of underdetermined linear systems are all manifestations of the same underlying algebraic property. We will call this algebraic property \textit{kernel atomicity}. It arises principally because homomorphic maps induce partitions of their domain space into cosets, `atoms' whose cardinalities are equal to that of the kernel. 


\end{abstract}\hspace{-1.5mm}\indent\indent {{\fontsize{7}{7}\selectfont {\bf Keywords:} Kernel, Orbit-Stabilzer Thm., First Isomorphism Thm., Homomorphism}}
\begin{center}
\line(1,0){300}
\end{center}

\newpage

\section{Introduction: Kernel Atomicity}

The kernel of an algebraic map is the set of elements that are mapped onto the identity. When this map is a homomorphism, it induces a partition of the domain space into cosets, `atoms' whose cardinalities are equal to the cardinality of the kernel. Here, we refer to this as kernel atomicity. It derives principally because composition of any given element of the domain by any element in the kernel has no effect on the given element's image. 

In the reminder of the paper we proceed as follows: In Subsection~(\ref{first-isomorphism}) we state the first isomorphism theorem; in Subsection~(\ref{orbit-stabilizer}) we state the orbit stabilizer theorem; in Section~(\ref{manifestations}) we present the manifestations of kernel atomicity, specifically, in Subsection~(\ref{manif:underdetermined}) we show how the non-uniqueness of solutions of underdetermined linear systems is a manifestation of kernel atomicity, in Subsection~(\ref{manif:first isomorphism}) we show how the first isomorphism is a manifestation of kernel atomicity, and in Subsection~(\ref{manif: orbit}) we show how the orbit stabilizer theorem is a manifestation of kernel atomicity. In Section~(\ref{summary}) we conclude with a summary of the paper. 

\subsection{The First Isomorphism Theorem}
\label{first-isomorphism}
The isomorphism theorems have played a unifying role in algebra since their discovery by Emmy Noether and Richard Dedekind~(\cite{no1927,vaba1931}). The first isomorphism theorem states that for $G$ a group, $H$ a subgroup of $G$, $\pi$ a homomorphism from $G$ onto $H$, and $Ker(\pi)$ the kernel of $\pi$; the quotient group $G/ Ker(\pi)$ is isomorphic to $H$. 

\subsection{The Orbit Stabilizer Theorem}
\label{orbit-stabilizer}

Given a group $G$ and a set $X$, consider the group action $\Pi:G \longmapsto Bij(X,X)$, where $Bij(X,X)$ is the group of bijections of $X$ onto itself. $\Pi(g)=\pi(g)$, where $\pi(g):X \longmapsto X$. The stabilizer of $x$ in $G$, denoted $Stab_G(x)$, is the subgroup of G for which $x$ is a fixed point under the associated bijections. Symbolically, $Stab_G(x) = \{g\in G~|~ \pi(g)x = x\}$. The orbit of $x$ in $G$, denoted $Orb_G(x)$, is the set of image points of $x$ under the action of elements of $G$. Symbolically, $Orb_G(x)=\{y\in X~|~ y=\pi(g)x ~\mbox{for some}~ g \in G\}$. The orbit-stabilizer theorem states that the cardinality of the orbit of $x$ in $G$ is the cardinality of $G$ divided by the cardinality of $Stab_G(x)$. Symbolically, $|G|/|Stab_G(x)| = |Orb_G(x)|$.

\section{Manifestations of Kernel Atomicity}
\label{manifestations}
\subsection{Underdetermined Linear System}
\label{manif:underdetermined}

Consider a linear transformation $L$ which maps a non-trivial element, $x$, to the identity. i.e. $Lx=0$. Then the kernel of $L$ consists of scalar multiples of $x$. i.e. $Ker(L) = \{cx~|~c\in \Re \}$. Therefore $|Ker(L)|=\infty$. If we then consider the non-homogenous system of equations $Ly=b$ with solution $y \neq 0$. It follows that $y+cx$ is also a solution $\forall c \in \Re$, because $L(y+cx)= L(y) + L(cx) = L(y) + cL(x) = b +c*0 = b$. Hence the solution is not unique because the pull-back of every point in the image space has infinite cardinality. This infinite cardinality is directly inherited from the kernel.

\subsection{First Isomorphism Theorem}
\label{manif:first isomorphism}
Given the homomorphic onto map $\pi: G \longmapsto H$. $\pi(x) = \pi(xy)~\forall y \in Ker(\pi)$. This holds because $\pi(xy)=\pi(x)\pi(y)=\pi(x)e_H=\pi(x)$, where $e_H$ is the identity in $H$. Therefore the pull-back of every image point $\pi(x)$ in $H$ is at least as large as $Ker(\pi)$, and contains elements of the form $xKer(\pi)$. To establish that there are no other types of elements in the pull-back of $\pi(x)$, and to establish that the size of the pull-back is exactly equal to the cardinality of $Ker(\pi)$, consider a $y\neq x$ such that $\pi(x)=\pi(y)$. Then $\pi(x)[\pi(y)]^{-1}=e_H$, therefore $\pi(xy^{-1}) = e_H$, therefore $xy^{-1} \in Ker (\pi)$, therefore $x \in yKer(\pi)$. Similarly, $y \in xKer(\pi)$. $x$ and $y$ are therefore in the same coset by $Ker(\pi)$. This establishes the atomicity of the kernel, and consequently, an isomorphism between the quotient space $G/Ker(\pi)$ and $H$. 

\subsection{$\mathbf{Stab_G(x)}$ is a Kernel}
\label{manif: orbit}
Here we prove the orbit stabilizer theorem and show how it derives from kernel atomicity. Consider the homomorphism $\Pi:G \longmapsto Bij(X,X)$. Then $\forall$ bijection $\pi(g)$ such that $\pi(g)x=x$, the action of $\pi(g)$ is the same as that of the identity map restricted to $x$. From this x-centric view, i.e. restriction of the identity map to $x$, the resulting kernel is: $\{g \in G ~|~ \pi(g)x=x\} = Stab_G(x)$. Hence under restriction to $x$, the kernel of $\Pi$ is the stabilizer of $x$ in $G$. Symbolically, $Ker(\Pi|_x) \simeq Stab_G(x)$. For each element $\pi|_x(g)$ of $Bij(X,X)|_x$ which is not analogous to the identity map, $\pi(g)x\neq x$. Let $\pi(g)x=y$. Then $\forall h \in Stab_G(x)$, $\pi(g)\pi(h)x=\pi(g)x=y$. Therefore $\pi(gh)=y~ \forall h \in Stab_G(x)$. Therefore the pull-back of every image point $\pi(g)x$ is at least as large as $Stab_G(x)$ and contains elements of the form $gStab_G(x)$. To establish that there are no other types of elements in the pull-back, and to establish that the size of the pull-back is exactly equal to the cardinality of $Stab_G(x)$, consider the following: $g,q \in G$ such that $g,q \notin Stab_G(x)$, and such that $\pi(g)x=\pi(q)x$. Therefore $[\pi(q)]^{-1}\pi(g)x=x$. Therefore $\pi(q^{-1}g)x=x$. Therefore $q^{-1}g \in Stab_G(x)$, hence $g \in qStab_G(x)$, and similarly $q \in gStab_G(x)$. $g$ and $q$ are therefore in the same coset by $Stab_G(x)$. This proves the orbit-stabilizer theorem that the cardinality of $Stab_G(x)$ must divide the cardinality of $Orb_G(x)$. And since $Ker(\Pi|_x) \simeq Stab_G(x)$, this is a manifestation of kernel atomicity.

\subsection{Homomorphism is Injective $\mathbf{\Leftrightarrow}$ Kernel $\mathbf{=\{e_G\}}$ }
\label{1-1}
Given a homomorphism $\pi$, the pull-back of every image point has a cardinality equal to that of the kernel. Therefore if the kernel contains only the identity, every image point pulls back to a unique element of the domain, and the map is necessarily injective. Conversely, if the map is injective then the pull-back of every image point, $\pi(x)$, is a singleton, $\{x\}$, for some $x\in G$. Since the pull-back of each image point is a coset in $G/Ker(\pi)$, it follows that $Ker(\pi)$ must be the identity element. Symbolically, $xKer(\pi)=\{x\}\Rightarrow Ker(\pi)=\{e_G\}$. 

\newpage
\section{Summary}
\label{summary}
On the surface, the first isomorphism theorem, the orbit stabilizer theorem, and the non-uniqueness of solutions of underdetermined linear systems appear  algebraically different. In this paper, we explicitly showed that they are manifestations of the same underlying algebraic property. This property which we termed kernel atomicity, refers to the partitioning of the domain of homomorphic maps into cosets, `atoms' whose cardinalities are equal to the cardinality of the kernel. The kernel of homomorphic maps thereby determines the granularity of partition, and consequently, the `atomic size'. As a corollary, when the kernel consists only of the identity element, the homorphism is necessarily injective.   

\newpage
\bibliographystyle{plainnat}
\bibliography{/home/odaibo/Dropbox/Research/Spatio_retinal_map/mybibliography}

\newpage
\section*{Author Biography}

{\footnotesize{
Dr. Stephen G. Odaibo is Chief Scientist and Founder of Quantum Lucid Research Laboratories, an Independent Research Institute dedicated to using a Computational Neuroscience approach to find a cure for blindness. Dr. Odaibo is also the 2014-2015 Medical Retina Fellow at the University of Michigan --Ann Arbor. He received the 2005 Barrie Hurwitz Award for Excellence in Clinical Neurology from Duke University School of Medicine where he topped the class in Neurology.

Dr. Odaibo is a Mathematician, Computer Scientist, Physicist, Neuroscientist, and Physician. He obtained a B.S. in Mathematics (UAB, 2001), M.S. in Mathematics (UAB, 2002), M.S. in Computer Science (Duke, 2009), and Doctor of Medicine (Duke, 2010). Dr. Odaibo completed his internship in Internal Medicine at Duke University Hospital and his residency in Ophthalmology at Howard University Hospital. He is author of the book, ``Quantum Mechanics and the MRI Machine'' (Symmetry Seed Books, Oct 2012).

He invented the Trajectron method with which he provided the first quantitative demonstration of non-paraxial light bending within the human cornea. His other awards and recognitions include: He won the 2013 Best Resident Research Presentation Award at the 23rd Annual Washington Retina Symposium; in 2012 he was selected as a Featured Alumnus of the Mathematics Department at UAB; and his cornea paper was selected by MIT Technology Review as one of the best papers from Physics or Computer science submitted to the arXiv the first week of Oct 2011.

Dr. Odaibo's research interests are at the fusion of Mathematics, Computer Science, Physics, and the Neuro-visual Sciences, with a special focus on the representation of motion in the mammalian visual cortex. His clinical interests are in the diagnoses and Medical/Laser management of retinal disease. He is also interested in developing the next generation of retinal and visual function imaging modalities.

He loves his wife, Lisa, his family and friends, and studying the bible.}}

\end{document}